\documentclass[a4paper, 12pt]{article}
\usepackage{hilbert}

\begin{document}

    \title{On Hilbert's Tenth Problem}
    \author{Michael Pfender\footnote{michael.pfender@alumni.tu-berlin.de}}
    \date{March 2014, last revised \today}
    \maketitle

\abstract{
Using an iterated Horner schema for evaluation of diophantine polynomials, 
we define a partial $\mu$-recursive ``decision'' algorithm \emph{decis}
as a ``race'' for a first \emph{nullstelle} versus a first (internal)
\emph{proof} of non-nullity for such a polynomial -- within a given theory $\T$
extending Peano Arithm\'etique $\PA$. If $\T$ is \emph{diophantine sound}, \ie
if (internal) \emph{provability} implies \emph{truth} -- for diophantine
formulae --, then the $\T$-map \emph{decis} gives \emph{correct results}
when applied to the codes of polynomial inequalities 
$D(x_1, \dots ,x_m) \neq 0.$ The additional hypothesis that $\T$ be
\emph{diophantine complete} (in the syntactical sense) would guarantee 
in addition termination of \emph{decis} on these formula, \ie \emph{decis}
would constitute a \emph{decision algorithm} for diophantine formulae in
the sense of Hilbert's 10th problem. From Matiyasevich's impossibility for
such a decision it follows, that a consistent theory $\T$ extending $\PA$
cannot be both diophantine sound and diophantine complete. We infer from
this the existence of a 
\emph{diophantine} formulae which is
undecidable by $\T$.
Diophantine correctness is inherited by the \emph{diophantine completion}
$\tildeT$ of $\T,$ and within this extension \emph{decis} terminates on all 
externally given diophantine polynomials, correctly. Matiyasevich's theorem 
-- for the strengthening $\tildeT$ of $\T$ -- then shows that $\tildeT,$ and 
hence $\T,$ cannot be diophantine sound. 
But since the internal consistency formula
$\mathrm{Con}_\T$ for $\T$ implies -- within $\PA$ -- diophantine soundness of $\T,$
we get $\PA \derives \neg \mathrm{Con}_\T,$ in particular $\PA$ must derive its own
internal inconsistency formula.
}


\section*{Overview}

\begin{enumerate}
 \item Consider a theory $\T$ with quantifiers and having terms for all 
  primitive recursive 
  maps (``\pr maps''); 
  so $\T$ is to be \NAME{Peano} Arithm\'etique $\PA$ or one of 
  $\PA$'s extensions, \eg $\mathbf{ZF}$ or $\mathbf{NGB}.$ 

 \item Obtain the theory $\tildeT$ by adding to $\T$ the axiom $\neg\mathrm{Con}_{\T}$
  of \emph{internal inconsistency}. By \NAME{G\"odel}'s second incompleteness 
  theorem, $\tildeT$ is consistent relative to $\T.$
        
 \item $\T$ admits a $\mu$-recursive, partially defined ``algorithm'' 
  $\decis$ aimed at deciding \emph{$\T$-internal} (\NAME{G\"odel}
  numbers of) \pr predicates.

 \item By \emph{internal semantical completeness} of $\tildeT$ with respect to 
  \pr predicates
  -- involving \emph{evaluation} of (\NAME{G\"odel} numbers of) internal 
    \pr predicates
     -- it is shown that in $\tildeT$ the
  \emph{partial} $\mu$-recursive $\T$-map $\decis$ is in fact \emph{total,} 
  and that it
  gives \emph{correct} results -- the latter for arguments $p$ of form
  $p = \code{\varphi},$ $\varphi = \varphi(n)$ a \pr predicate, 
  $\code{\varphi} \in \N$ its internal \NAME{G\"odel} number.

 \item within $\tildeT,$ $\decis$ decides in particular (systems of) 
  \emph{diophantine equations.}
        
 \item \NAME{Matiyasevich}'s negative result concerning this decision
  problem of \NAME{Hilbert} is a theorem of $\T,$ a fortiori of $\tildeT.$

 \item This contradiction shows $\tildeT,$ hence also $\T,$ to be inconsistent:
  ``unbounded formal
    quantification is incompatible with infinity.''
\end{enumerate}

\section{Decision}

Crucial for the present approach to \NAME{Hilbert}'s decision problem
is availability -- within $\T$ -- of a ($\mu$-recursive) 
\emph{evaluation} map
 $\ev: \N \times \N \supset |\N,2|_{\PR} \times \N \to 2$
on the $\T$-internal (primitive recursively decidable) \emph{set} 
$|\N,2|_{\PR} \subset \N$ of 
\NAME{G\"odel} numbers (``codes'') 
of \pr predicates.
(Primitive recursive \emph{predicates} are viewed as \pr \emph{map terms} 
with codomain 
$2 \subset \N).$ This evaluation map 
$ev$ is defined in $\T$ by 
\emph{(nested) double recursion} \`a la 
\NAME{Ackermann}, see \NAME{P\'eter} 1967, and 
satisfies the characteristic equation 
  $$\ev(\code{\varphi}, n) = \varphi(n)$$
for \pr predicates $\varphi = \varphi(n)$ of $\T,$ \cf Appendix. Here
$\code{\varphi} \in |\N,2|_{\PR} \subset \N$ is $\varphi$'s $\T$-internal
\NAME{G\"odel} number.
 
\emph{Define} now the partial $\mu$-recursive ``decision'' $\tildeT$-map
 $$\decis = \decis(p): |\N,2|_{\PR} \partialmapto 2$$
hoped for deciding (internal) \pr predicates $p,$ \ie 
 $p \in |\N,2|_{\PR} \subset \mathit{formulae}_{\tildeT} = 
  \mathit{formulae}_{\T} \subset \N,$
    via the two ``antagonistic'' \emph{termination indices}
 $$\mu_{\ex}(p), \mu_{\thm}(p):
      |\N,2|_{\PR} \to \N \union \{\infty\} \text{ as follows:}$$
    \begin{align*}
      \mu_{\mathit{ex}}(p) & := \mu \{ n: \ev(p, n) = 0 \} \quad 
             \text{``minimal counter\emph{ex}ample''} \\
 & =_{\mathit{def}} \begin{cases}
           \min \{ n: \ev(p,n) = 0 \} 
                \quad \text{if} \quad \exists n (\ev(p,n) = 0) \\
           \infty \ (\corr undefined) 
                \quad \text{if} \quad \forall n (\ev(p,n) = 1);
      \end{cases}
    \end{align*}

the \emph{theorem index} 
$\mu_{\thm}(p) \in \N \union \{\infty\}$ of
$p \in |\N,2|_{\PR}$ is defined by
 $$\mu_\thm(p) := \mu\{ k: \thm_{\tildeT} (k) = p \};$$
here the \pr enumeration 
$\thm_{\tildeT} = \thm_{\tildeT}(k): \N \to  
\mathit{formulae}_{\T} \subset \N$
is the \emph{$\tildeT$-internal} version of the
\emph{metamathematical} enumeration of all (\NAME{G\"odel} numbers of) 
$\tildeT$-theorems; enumeration is
\emph{lexicographic} by ``length of shortest proof''.
  
Finally, we define the -- a priori partial -- \emph{$\mu$-recursive} $\T$-map 
  $$\decis = \decis(p): |\N,2|_{\PR} \partialmapto 2 \text{ by }$$
  \begin{displaymath}
    \decis(p) =
    \begin{cases}
      0 \text{ if } \mu_\ex(p) < \infty  
           \quad \text{(``counter\emph{ex}ample'')} \\
      1 \text{ if } \mu_\ex(p) = \infty \text{ and } 
                    \mu_\thm(p) < \infty \quad \text{(``theorem'')} \\
      \infty \text{ otherwise, \ie if } 
            \mu_{\thm}(p) = \mu_{\mathit{ex}}(p) = \infty. 
    \end{cases}
  \end{displaymath}


For proving $\decis$ to be totally defined within 
$\tildeT$ = $\T+\neg\mathrm{Con}_{\T}$ 
we rely on the following

\bigskip

\textbf{Lemma} (Internal Semantical Completeness): 
  $$\tildeT \vdash \forall n (ev(p,n) = 1) \implies 
        \exists k (\thm_{\tildeT}(k) = p)$$
with $p$ free on $|\N,2|_{\PR},$ in closed form:
   $$\tildeT \vdash (\forall p \in |\N,2|_{\PR}) [\forall n (ev(p,n) = 1) 
                                 \implies \exists k (\thm_{\tildeT}(k) = p)].$$

\textbf{Proof:} One of the equivalent $\T$-formulae expressing internal
inconsistency of $\T$ is
 $$\neg\mathrm{Con}_{\T} = (\forall f \in \mathit{formulae}_{\T})
   (\exists k) (\thm_{\T}(k) = f):$$
``every internal formula (its \NAME{G\"odel} number in $\T$) is \emph{provable}''
(emphasis from \NAME{G\"odel}). This gives in particular
 $$\tildeT \vdash \exists k (\thm_{\tildeT}(k) = p),$$
$p$ free on $|\N,2|_{\PR} \subset \mathit{formulae}_{\T} \subset \N,$ 
and hence -- trivially -- the assertion of the Lemma.

\bigskip

\textbf{Decision Lemma:}
    \begin{enumerate}
     \item within $\tildeT = \T+\neg \mathrm{Con}_{\T}$, the (a priori partial) 
      $\mu$-recursive
          \emph{decision-``algorithm''}
       $$\decis(p): |\N,2|_{\PR} \partialmapto 2$$
      is in fact \emph{totally defined}, with other words it \emph{terminates} 
      on all 
      internal \NAME{G\"odel} numbers $p \in |\N,2|_{\PR}.$

     \item For $\varphi = \varphi(n)$ a \pr predicate, 
      $\code{\varphi} \in |\N,2|_{\PR} \subset \N$ its $\T$-internal 
      \NAME{G\"odel} number,
      $\decis(\code{\varphi})$ gives -- in $\tildeT$ -- the \emph{correct} 
      result:
      \begin{itemize}
       \item $\tildeT \vdash \decis(\code{\varphi})=0 \iff \exists n
            (\neg \varphi (n)),$
       \item $\tildeT \vdash \decis(\code{\varphi})=1 \implies 
             \forall n (\varphi (n)).$
      \end{itemize}
    \end{enumerate}

    \textbf{Proof} of (i):
     \begin{align*}
          \tildeT \derives  [ \quad 
                        & \mu_\ex(p) = \infty  \\
                        & \iff \forall n (ev(p,n) = 1)  \\
                        & \implies \exists k (\thm_{\tildeT}(k) = p) \\
                        &  \quad\quad\quad \text{by internal semantical 
                                         completeness of $\tildeT$ above} \\
                        & \iff \mu_{\thm}(p) < \infty \quad ].
     \end{align*}
     Hence not both of $\mu_{\ex}(p), \mu_{\thm}(p)$ can be undefined. This 
     shows \emph{termination}
       $$\decis(p) \in \{0, 1\}$$
      of $\decis$ within $\tildeT$ for all (internal) \pr\ predicates $p$ 
      (\NAME{G\"odel} numbers thereof).
  
      Proof of (ii):
      \begin{align*}
          \tildeT \derives [ \quad
               & \decis(\code{\varphi}) = 0  \\
               & \iff \mu_\ex(\code{\varphi}) < \infty \\
               & \iff \exists n (\ev(\code{\varphi},n) = 0) \\
               & \iff \exists n (\varphi(n) = 0) 
                  \quad \text{by $\ev$'s evaluation property} \\ 
               & \iff \exists n (\neg \varphi(n)) \quad ] 
                               \text{ as well as } \\
          \tildeT \derives [ \quad
               & \decis(\code{\varphi}) = 1  \\
               & \implies \mu_\ex(\code{\varphi}) = \infty \\
               & \iff \forall n (\ev(\code{\varphi},n) = 1) \\
               & \iff \forall n (\varphi(n)) \quad ] \quad \text{\qed}
      \end{align*}

\section{Hilbert's 10th Problem revisited}


  \bigskip
  
  A system
  \begin{displaymath}
      D: \qquad 
      \begin{array}{ccc}
          D_{1}^{L}(x_1, \ldots, x_m) & = & D_{1}^{R}(x_1, \ldots, x_m) \\
          \vdots & & \vdots  \\
          D_{k}^{L}(x_1, \ldots, x_m) & = & D_{k}^{R}(x_1, \ldots, x_m)
      \end{array}
  \end{displaymath}
  of $k$ \emph{diophantine equations} -- see 
  \NAME{Matiyasevich} 1993, 1.1, 1.2, and 1.3 -- gives rise to a 
  \emph{\pr predicate}
  \begin{align*}
      & \varphi = \varphi(x_1, \ldots, x_m): \N^{m} \to 2 \text{ defined by} \\
      & \varphi(x_1, \ldots, x_m) =
        [D_{1}^{L} \neq D_{1}^{R} \vee \ldots \vee D_{k}^{L} \neq D_{k}^{R}]: 
        \N^m \to 2
  \end{align*}
having the property that $(x_1, \dots, x_m) \in \N^{m}$ 
is a solution to system $(D)$ iff it is a \emph{counterexample to} $\varphi,$
and $(D)$ has \emph{no} \emph{solution} (in natural numbers) iff $\varphi$ 
\emph{holds} for 
$(x_1, \dots, x_m) \text{ \emph{free} in } \N^{m}.$ 

\NAME{Cantor}'s p.r. enumeration 
$cantor_{m}: \N \to \N^{m}$
having a p.r. inverse
$cantor_{m}^{-1}: \N^{m} \to \N,$
$$\psi = \psi(n) := \varphi(\cantor_{m}(n)): \N \to 2$$
is a p.r. predicate of $\T$ such that 
$(x_1, \dots, x_m) \in \N^{m}$ 
\emph{solves} $(D)$ iff
$\cantor_{m}^{-1}(x_1, \dots, x_m) \in \N$ is a \emph{counterexample} to $\psi,$ and
$(D)$ is \emph{unsolvable} iff $\psi(n)$ \emph{holds} for 
$n \text{ \emph{free} in } \N.$
So from the Decision Lemma (for p.r. predicates) above we obtain:

\bigskip

\textbf{Decision Theorem:}
\begin{enumerate}
\item 
 Within the -- somewhat strange -- theory $\tildeT = \T+\neg\mathrm{Con}_{T},$ 
 the (partial) $\mu$-recursive map (the ``\emph{algorithm}'')
 $\decis: |\N,2|_{\PR} \partialmapto 2$ 
 \emph{decides} all (internal) primitive recursive predicates, in
 particular all (internal, a fortiori external) \NAME{G\"odel} numbers coding
 ``diophantine'' predicates as considered above, and hence decides internal, 
 a fortiori
 external (systems of) \emph{Diophantine equations.}

\item
 Since $\mu$-recursion and \NAME{Turing}-machines have equal 
 \emph{computation power} -- by the verified part of \NAME{Church}'s thesis -- 
 this means: Within $\tildeT,$ $\decis$ gives rise to a \NAME{Turing}
 machine $TM$ deciding all internally given as well as all externally given 
 Diophantine equations, \ie $\tildeT$ admits a \emph{positive} solution to 
 \NAME{Hilbert}'s 10th 
 problem.

\item
 On the other hand, \NAME{Matiyasevich}'s \emph{negative} 
 solution to this problem 
 needs as a formal framework $\T$ just Arithm\'etique $+ \exists.$  

\item The latter two results -- \NAME{Matiyasevich}'s \emph{negative}
 $\T$-\emph{theorem} and our \emph{positive} $\tildeT$-\emph{theorem} 
 contradict each other in the stronger theory
 $\tildeT$. This shows $\tildeT$ to be \emph{inconsistent}. 

\item
 \NAME{G\"odel}'s consistency of $\neg\mathrm{Con}_{\T}$ relative to $\T$ then entails 
 inconsistency of $\T,$ whence in particular inconsistency of 
 \NAME{Peano} Arithm\'etique $\PA$ and of the classical set theories.

\end{enumerate}

\bigskip

\textbf{Corollary:} Since \NAME{Matiyasevich} 1993 makes essential use of
formal (existential) quantification for ``unsolving'' \NAME{Hilbert}'s 10th 
problem, this only decision problem on \NAME{Hilbert}'s list is again
open -- for treatment within the framework of a suitable \emph{constructive}
foundation for Arithmetic.




\section{Appendix: Evaluation}


In section 2 we made appeal to availability in $\T$ of an \emph{evaluation} 
$\ev = \ev(p,n)$ 
of (internal) \pr predicate codes $p$ satisfying 
  $$\ev(\code{\varphi},n) = \varphi(n)$$
for (``external'') \pr predicates $\varphi: \N \to 2$ in $\T.$
We identify a \pr
predicate $\varphi = \varphi(n)$ of $\T$ with
its associated \pr map term $\varphi = \varphi(n): \N \to 2,$
since we want to define the evaluation of (internal) \pr predicates by 
restriction of an evaluation of \emph{all} internal \pr map terms out of the set
$|\N,2|_{\PR} \subset \N$ of (internal) \pr map terms from $\N$ to $2.$

For \emph{defining} this map term evaluation $\ev$ by 
\emph{(nested) double recursion} \`a la 
\NAME{Ackermann} (\cf \NAME{P\'eter} 1967) we need a 
\emph{universal set} (object)
  $$\U = \N^{\,(*)}$$
of all \emph{nested pairs} of natural numbers, and hence containing all
$\PR$-objects $1,\N,\dots, A,\dots, B, A \times B, \dots$ as \emph{disjoint}
(exception: $1 \subset \N$) \pr decidable subsets.

This set $\N^{\,(*)}$ is directly available in \emph{set theory.} Within
\NAME{Peano} Arithm\'etique, it can be ``constructed'' via \emph{coding}
as a decidable subset of $\N.$

\bigskip

\textbf{Definition}: \emph{Evaluation}
 $$\ev = ev(u,a): \N \times \N^{\,(*)} \supset 
                     \mathit{PR} \times \N^{\,(*)} \to \N^{\,(*)} $$
of the \emph{internal} (\NAME{G\"odel} numbers of) \pr \emph{maps} 
$u,v,w \in \mathit{PR} \subset \N,$  
on binary nested tupels
$a,b,c \in \N^{\,(*)}$ of natural
numbers is now defined by \emph{(nested) double recursion} with principal 
recursion parameter ``operator-depth'' $\depth(u)$ of $u$ as follows:

\begin{itemize}

\item basic internal map terms $\code{0}, \code{s}, \code{id}, \code{!},
       \code{\Delta}, \code{\Theta}, \code{\ell}:$
\spiegel \quad $\ev(\code{0},0) = 0 = 0(0) \in \N$ ``zero map'',

\spiegel \quad $\ev(\code{s},n) = n+1 = s(n) \in \N$ ``successor map'',

\spiegel \quad $\ev(\code{\id},a) = a = id(a)$ ``identity'',

\spiegel \quad $\ev(\code{!},a) = 0 = \,!(a) \in 1 \subset \N$ ``terminal map'',

\spiegel \quad $\ev(\code{\Delta},a) = (a,a) = \Delta(a)$ ``diagonal'',

\spiegel \quad $\ev(\code{\Theta},(a,b)) = (b,a) = \Theta(a,b)$ 
                  ``transposition'',

\spiegel \quad $\ev(\code{\ell},(a,b)) = a = \ell(a,b)$ ``left projection''.


This defines $\ev$ on $\PR$'s (map-)\emph{constants}, $\depth$ of these
``basic'' map terms is set to 1.

We now define $\ev$ on compound internal \pr map terms:

\item internally \emph{composed} $v \code{\circ} u$:
 $$\ev(v \code{\circ} u, a) = \ev(v,\ev(u,a)).$$
This definition is \emph{legitimate}, since
\begin{align*}
 & \depth(u), \depth(v) < \depth(v \code{\circ} u) \\ 
 & =_{\mathit{def}} \max(\depth(u),\depth(v))+1 \in \N;
\end{align*}

\emph{Example:} 
\begin{align*}
 & \ev(\code{s} \code{\circ} \code{s} \code{\circ} \code{s},s(0)) \\
 & = \ev(\code{s},\ev(\code{s},\ev(\code{s},s(0)))) \\
 & = ((s(0)+1)+1)+1 = 4.
\end{align*}



\item cylindrified $\code{\id} \code{\times} v:$
 $$\ev(\code{id} \code{\times} v, (a,b)) = (a,\ev(v,b)),$$
``evaluation in the second component''. 

\emph{legitimacy of this definition:} 
 $$depth(v) < depth(\code{id} \code{\times} v) =_{\mathit{def}} \depth(v)+1.$$

\item internally \emph{iterated} $u^{\S}$:
  \begin{align*}
   \ev(u^{\S},(a,0)) & = a, \\
   \ev(u^{\S},(a,n+1)) 
    & = \ev(u,\ev(u^{\S},(a,n))).
  \end{align*}
\end{itemize}

This last case is in fact a \emph{(nested) double recursion} 
\`a la \NAME{Ackermann},
since the \emph{internally iterated} $u^{\S}$ of $u$ is evaluated 
in a \pr manner with respect
to the second parameter $n \in \N$ -- which is to count the iteration loops 
still to be performed. The principal recursion parameter is (internal)
operator-depth 
$\depth = \depth(u): \N \supset \mathit{PR} \to \N,$ 
in particular in this last case  
$\depth(u^{\S}) =_{\mathit{def}} \depth(u)+1.$

Each primitive recursive map can be generated from the basic maps
$0,s,id,!,\Delta,\Theta,\text{ and } \ell$ by composition, cylindrification
and iteration: \emph{substitution} is realized via composition with the
induced $(f,g) = (f,g)(c) = (f(c),g(c))$ which in turn is obtained via
diagonal, cylindrification, transposition, and composition. 
Since iteration $g^\S$ then gives the (``full'') schema of 
primitive recursion (see \NAME{Freyd 1972}, \NAME{Pfender} et al. 1994), 
$\ev$ in fact evaluates
all \NAME{G\"odel} numbers of (internal) \pr map terms, recursively given
in the above way. 

Let us call $\PR+ev$ the extension of $\PR$ by a (formal) map
 $$\ev = \ev(u,a): 
   \N \times \N^{\,(*)} \supset \mathit{PR} \times \N^{\,(*)} \to  \N^{\,(*)}$$
satisfying the above 2-recursive system for $\ev.$

For our ``set'' theory $\T$ we now prove the following

\bigskip

\textbf{Evaluation Lemma:}
For primitive recursive $f: \N^{\,(*)} \supset A \to B \subset \N^{\,(*)}$ 
in $\T,$ 
$\T$ extending $\PR+\ev,$ we have
 $$\ev(\code{f},a) = f(a): A \to B,$$
in particular for $\varphi: \N \to 2$ 
   (the map term representing) a \pr predicate of $\T:$
 $$\ev(\code{\varphi},n) = \varphi(n): \N \to 2, 
       \quad n \text{ free variable on } \N.$$

\textbf{Proof} by \emph{external} (``metamathematical'') \emph{induction} 
on the operator-depth 
$\bfdepth(f) \boldsymbol{\in} \boldsymbol{N}$
of $f$ varying on $\PR \boldsymbol{\subset} \boldsymbol{N},$ 
in case of an iterated
$f = g^\S(a,n): A \times \N \to A$ 
this external induction will be combined with an internal induction
on the \emph{iteration parameter} $n \in \N.$ 
$\bfdepth: \PR \boldsymbol{\to} \boldsymbol{N}$
is the external primitive recursive ``twin'' of $\depth: \mathit{PR} \to \N$ 
above;
it is characterised by 
$\depth(\code{f}) = \num(\bfdepth(f))$
for $f: A \to B \text{ in } \PR \boldsymbol{\subset} \T.$ Here 
$\num = \num(\boldsymbol{n})\boldsymbol{:} \boldsymbol{N} \boldsymbol{\to}
 \T(1,\N)$ maps each \emph{external} natural number $\boldsymbol{n}$ 
into its corresponding
$\T$-\emph{numeral}, as defined \eg in set theory by associating 
\NAME{von Neumann} numerals.

\begin{itemize}

\item Anchoring: the assertion holds for the \emph{basic} maps 
  $0,\dots, \ell$ (with $\bfdepth$ set to 
  $\boldsymbol{1} \boldsymbol{\in} \boldsymbol{N}$) 
  just by definition of $\ev.$ 

\item composition case $f = h \circ g: A \to B \to C:$
  \begin{align*}
    & \ev(\code{f},a) = \ev(\code{h \circ g},a)\\
    & = \ev(\code{h} \code{\circ} \code{g},a)
         \quad \text{since } \code{\circ} 
         \text{ internalizes } \lq \circ \rq \\
    & = \ev(\code{h},\ev(\code{g},a))
         \quad \text{by definition of } \ev \\
    & = \ev(\code{h},g(a)) \quad \text{by recursion hypothesis on } g\\ 
    & \quad \quad \text{since }
       \bfdepth\boldsymbol{(} g \boldsymbol{)} 
       \boldsymbol{<} \bfdepth\boldsymbol{(} f \boldsymbol{)} \\
    & = h(g(a)) \quad \text{by recursion hypothesis on } h \\ 
    & \quad \quad \text{since } 
       \bfdepth\boldsymbol{(} h \boldsymbol{)} 
       \boldsymbol{<} \bfdepth\boldsymbol{(} f \boldsymbol{)} \\
    & = (h \circ g)(a) = f(a).
  \end{align*}

\item case $f = id \times g: A \times B \to A \times C$ a cylindrified map:
  \begin{align*}
    & \ev(\code{f},(a,b)) = \ev(\code{id \times g},(a,b)) \\
    & = \ev(\code{id} \code{\times} \code{g},(a,b)) \\ 
      & \quad \quad 
          \text{since } \code{\times} \text{ is to internalize } \times \\
    & = (a,\ev(\code{g},b)) \quad \text{by definition of } \ev \\
    & = (a,g(b)) \quad \text{by recursion hypothesis on } g \\
    & \quad \quad \text{since } 
      \bfdepth(g) 
       \boldsymbol{<} \bfdepth(f) \\
    & = (id \times g)(a,b) = f(a,b).
  \end{align*}

\item The remaining case -- not quite so simple -- 
    is that of an \emph{iterated}
   $f = g^\S: A \times \N \to A$ of a (\pr) endo map $g: A \to A,$
   $g^\S$ characterized by 
    $$g^\S(a,0) = a, \quad g^\S(a,n+1) = g(g^\S(a,n)):$$
   the assertion of the Lemma holds in this last case too, 
   since -- ``anchoring'' $n = 0$ for internal induction:
   \begin{align*}
         & \ev(\code{f},(a,0)) = \ev(\code{g^\S},(a,0)) \\ 
         & = \ev(\code{g}^{\S},(a,0)) = a 
                 \quad \text{since } \langle \_ \rangle^{\S} 
                 \text{ internalizes } \myHat{(\_)} \\
         & = g^\S(a,0) = f(a,0)
   \end{align*}
   -- as well as (internal induction step, using the external 
   recursion hypothesis):
   \begin{align*}
     & \ev(\code{f},(a,n+1)) = \ev(\code{g^\S},(a,n+1)) \\
     & = \ev(\code{g}^{\S},(a,n+1)) 
              \quad \text{since } \langle \_ \rangle^{\S}
              \text{ internalizes } \myHat{(\_)} \\
     & = \ev(\code{g},\ev(\code{g}^{\S},(a,n))) \\
       & \quad \quad \text{by (internal) inductive definition of } \ev \\
       & \quad \quad \text{in the present case } v = u^{\S} = \code{g}^{\S} \\
     & = \ev(\code{g},\ev(\code{g^\S},(a,n)))
       \quad \text{by $\langle \_ \rangle^{\S}$ internalizing } \myHat{(\_)} \\
     & = \ev(\code{g},g^\S(a,n))
         \quad \text{by (internal) induction hypothesis on } n \\
     & = g(g^\S(a,n))
         \quad \text{by (external) recursion hypothesis on } g \\
     & \quad \quad \text{since } \bfdepth(g)
         \boldsymbol{<} \bfdepth(f)) \\
     & = g^\S(a,n+1) = f(a,n+1) \text{ by definition of the
       iterated } g^\S \text{ \qed}
   \end{align*} 
\end{itemize}

\bigskip \bigskip
{ \setlength{\parindent}{0em} \addtolength{\parskip}{1ex}

  \begin{large}
    References
  \end{large}
  

  
  \NAME{S.\ Eilenberg, C.\ C.\ Elgot} 1970: \emph{Recursiveness}.
  Academic Press.
  
  
  
  
  \NAME{P.\ J.\ Freyd} 1972: Aspects of Topoi. \emph{Bull.\ Australian 
   Math.\ Soc.\ } 7, 1-76.
  
  \NAME{K.\ G\"odel} 1931: \"Uber formal unentscheidbare S\"atze der
  Principia Mathematica und verwandter Systeme I. \emph{Monatsh.\ der
    Mathematik und Physik} 38, 173-198.
  
  \NAME{R.\ L.\ Goodstein} 1971: \emph{Development of Mathematical
    Logic}, ch. 7: Free-Variable Arithmetics.  Logos Press.
  
  \NAME{D.\ Hilbert} 1900: Mathematische Probleme. Vortrag. Quoted in
  \NAME{Matiyasevich} 1993.
  
  
  
  
  
  
  \NAME{F.\ W.\ Lawvere} 1964: An Elementary Theory of the Category of
  Sets.  \emph{Proc.\ Nat.\ Acad.\ Sc.\ USA} 51, 1506-1510.
  


  
  \NAME{Y.\ V.\ Matiyasevich} 1993: \emph{Hilbert's Tenth Problem}.
  The MIT Press.
  
  \NAME{R.\ P\'eter} 1967: \emph{Recursive Functions}. Academic Press.

  
  \NAME{M.\ Pfender, M.\ Kr\"oplin, D.\ Pape} 1994: Primitive
  recursion, equality, and a universal set. \emph{Math.\ Struct.\ in
    Comp.\ Sc.\ } 4, 295-313.
  
  
  

  
  
  
  \NAME{C.\ Smorynski} 1977: The incompleteness theorems. Part D.1,
  pp. 821-865 in \NAME{J.\ Barwise} ed. 1977: \emph{Handbook of
    Mathematical Logic}.  North Holland.



\end{document}

\bigskip

  Address of the author:\\
  \NAME{M. Pfender}\\
  Institut f\"ur Mathematik\\
  Technische Universit\"at Berlin MA 1-1\\
  Str. d. 17.Juni 136\\
  D-10623 Berlin\\
  
  pfender@math.tu-berlin.de
}

\end{document}